\documentclass[10pt]{article}
\usepackage{amssymb, amsmath, url, graphicx, setspace, geometry}
\usepackage{calrsfs}
\usepackage{wasysym, xcolor}
\usepackage{booktabs}
\usepackage{array}
\usepackage{enumitem}

\DeclareMathOperator{\os}{os}
\DeclareMathOperator{\lcm}{lcm}
\DeclareMathOperator{\expo}{exp}

\title{\bf On some problems regarding $LCM$-groups}
\author{Mihai-Silviu Lazorec}
\date{September 7, 2026}

\begin{document}

\maketitle

\begin{abstract}
Let $G$ be a finite group and denote by $o(g)$ the order of an element $g\in G$. We say that $G$ is an $LCM$-group if $o(x^ny)$ is a divisor of the least common multiple of $o(x^n)$ and $o(y)$ for all $x, y\in G$ and $n\in\mathbb{N}$. This paper extends some existing results on $LCM$-groups, such as the structure of a minimal non-$LCM$-group, and establishes criteria for $G$ to be an $LCM$-group or a nilpotent group. We also prove  that, in general, a minimum cover of a finite set of $LCM$-groups is not an $LCM$-group, and we answer two questions posed by M. Amiri.           
\end{abstract}



\section{Introduction}

All groups considered in this paper are finite. Let $G$ be a group and let $k\in\mathbb{N}$. The order of an element $x\in G$, the exponent of $G$, the set of prime divisors of $|G|$ and the Fitting subgroup of $G$ are denoted by $o(x)$, $\expo(G)$, $\pi(G)$ and $Fit(G)$, respectively. If $G$ is a $p$-group, then $n_{p^i}(G)$ is the number of elements of order $p^i$ in $G$ and $\Omega_k(G)=\langle x\in G\mid x^{p^k}=1 \rangle$ are the omega subgroups of $G$.  $C_k$ denotes the cyclic group of order $k\geq 2$. Other groups that appear in the paper are:
\begin{itemize}[itemsep=0pt, topsep=2pt]
\item the modular $p$-group $M_{p^k}=\langle a,b\mid a^{p^{k-1}}=b^p=1,\, ba=a^{p^{k-2}+1}b\rangle$ ($k\ge 4$ if $p=2$; $k\ge 3$ if $p$ odd);
\item the dihedral group $D_{2k}=\langle a,b\mid a^k=b^2=1,\, ba=a^{k-1}b\rangle$ ($k\ge 3$);
\item the generalized quaternion group $Q_{2^k}=\langle a,b\mid a^{2^{k-1}}=1,\, b^2=a^{2^{k-2}},\, ba=a^{2^{k-1}-1}b\rangle$ ($k\ge 3$);
\item the quasidihedral group $QD_{2^k}=\langle a,b\mid a^{2^{k-1}}=b^2=1,\, ba=a^{2^{k-2}-1}b\rangle$ ($k\ge 4$).
\end{itemize}

The study done on the sum of element orders of $G$, defined as
$$\psi(G)=\sum\limits_{x\in G}o(x),$$ is part of the more general research topic of  describing groups' properties via element orders. The literature on this invariant is  vast, so we only recall that there are several papers (such as \cite{1}, \cite{4}, \cite{5}, \cite{6}, \cite{10}, \cite{11}, \cite{18}, \cite{21}) where $\psi(G)$ was used to describe the nature (cyclicity, nilpotency, and so on) and structure of $G$. 

In this paper, we obtain some results mainly  referring to $LCM$-groups. These groups are defined by a divisibility relation written in terms of element orders, as shown  below. Hence, our work is connected with the same research topic of characterizing the properties of a group by its element orders. In the rest of this section, we recall some  concepts that are used or investigated in this paper and we briefly outline our objectives. 

We consider the sets
$$LCM(G)=\{ x\in G \mid o(x^ny)\mid \lcm(o(x^n), o(y)), \forall \ n\in\mathbb{N}, \forall \ y\in G \}$$
$$CP_2(G)=\{ x\in G \mid o(xy)\leq \max\{ o(x), o(y)\}, \forall \ y\in G\}.$$
We say that $G$ is a:
\begin{itemize}
\item $LCM$-group if $LCM(G)=G$;
\item $CP_2$-group if $CP_2(G)=G$.
\end{itemize}
If $G$ is a $p$-group, we also consider the sets $$P_2^k(G)=\{ x\in \Omega_k(G) \mid o(x)\leq p^k\},$$
and we say that $G$ is a $P_2^*$-group if $P_2^k(G)=\Omega_k(G),$ for all $k\in\mathbb{N}.$ We mention that these 3 properties defined via element orders were mainly investigated by Amiri and Lima in \cite{3}, T\u arn\u auceanu in \cite{20}, and Mann in \cite{15}, respectively. We also recall that Mann said that a $p$-group $G$ is a $P_2$-group if all sections of $G$ are $P_2^*$-groups.

If $\mathcal{F}$ is a finite set of groups, then $G$ is called an $\mathcal{F}$-cover if,   for each $F\in\mathcal{F}$, $G$ has a subgroup isomorphic to $F$. If $G$ has the smallest order among all $\mathcal{F}$-covers, then $G$ is called a minimum $\mathcal{F}$-cover. These concepts were recently studied by Cameron, Craven, Dorbidi, Harper and Sambale in \cite{8}. Any set $\mathcal{F}$ has a (minimum) $\mathcal{F}$-cover since the direct product of the groups in $\mathcal{F}$ is an $\mathcal{F}$-cover.

The order sequence of $G$, denoted by $\os(G)$, is a non-decreasing sequence formed of the element orders of $G$. Explicitly, if $G=\{ x_1, x_2, \ldots, x_k\}$, then the order sequence of $G$ is
$$\os(G)=(o(x_1), o(x_2), \ldots, o(x_k)), \text{ \ where \ } o(x_i)\leq o(x_{i+1}), \ \forall \ i \in \{ 1, 2, \ldots, k-1\}.$$
For the ease of writing, if the order sequence has a greater size and it includes duplicates, we write some pairs such as $(o(x), i)$ into $\os(G)$ meaning that $o(x)$ has a multiplicity of $i$ (i.e. $o(x)$ appears in $\os(G)$ for $i$ times). For instance, for the cyclic group $C_k$ of order $k$, we have
$$\os(C_k)=((1, 1), (d_1, \varphi(d_1)), (d_2, \varphi(d_2)), \ldots, (d_l, \varphi(d_l)), (k, \varphi(k))),$$ where $d_1, d_2, \ldots, d_l$ are the non-trivial divisors of $k$ and $\varphi$ is Euler's totient function. We refer the reader to papers \cite{9} and \cite{14}, by Cameron and Dey and by Lazorec, respectively, for more details, properties and results concerning order sequences.

Our main objectives are:
\begin{itemize}
\item describing the structure of minimal non-$LCM$-groups;
\item obtaining some criteria for $G$ to be an $LCM$-group or a nilpotent group by means of two ratios that are defined in terms of $LCM(S)$, where $S$ is a section of $G$;
\item showing that, in general, if $\mathcal{F}$ is formed of $LCM$-groups, none of its minimum $\mathcal{F}$-covers is an $LCM$-group; 
\item answering two questions that include hypotheses on the $LCM$ property, order sequences and sum of element orders; both questions were posed by Amiri in \cite{2}.
\end{itemize}
More details about these problems, results and proofs are given in the following sections.
 
\section{Minimal non-$LCM$-groups}

We recall that given a class $\mathcal{X}$ of groups, it is said that $G$ is a minimal non-$\mathcal{X}$ group if $G\not\in \mathcal{X}$ and $H\in\mathcal{X}$ for any proper subgroup $H$ of $G$. For instance, choosing $\mathcal{X}$ to be the class of nilpotent groups, one would work with minimal non-nilpotent groups (also known as Schmidt groups).

The first two results below list known facts about Schmidt groups (see \cite{16},  \cite{17}, Lemma 2.3 of \cite{13} and the proof of Theorem 1 in \cite{20}) and a  characterization of $LCM$-groups (see Theorem 2.6 of \cite{3}). The third lemma shows that the $LCM$ property is inherited by subgroups and direct products.\\

\textbf{Lemma 2.1.} \textit{Let $G$ be a Schmidt group. Then the following statements hold:
\begin{itemize}
\item[i)] $G=P\rtimes Q$ where $P$ is a normal Sylow $p$-subgroup, $Q=\langle y \rangle$ is a cyclic Sylow $q$-subgroup of $G$, $y^q\in Z(G)$ and $p, q$ are distinct primes;
\item[ii)] $Z(G)=\Phi(G)=\Phi(P)\times \langle y^q\rangle$;
\item[iii)] $\frac{P}{\Phi(P)}\cong C_p^r$, where $r$ is the multiplicative order of $p$ modulo $q$;
\item[iv)] $\frac{G}{Z(G)}$ is a Frobenius group such that the Frobenius kernel is elementary abelian and the Frobenius complement is of prime order; 
\item[v)] The maximal subgroups of $G$ are $P\times \langle y^q\rangle$ and $\Phi(P)Q_i$, where $Q_i$ is a conjugate of $Q$ and $i\in \{1, 2, \ldots, n_q\}$.
\end{itemize}}

\textbf{Lemma 2.2.} \textit{Let $G$ be a group. Then $G$ is an $LCM$-group if and only if $G$ is nilpotent and all its Sylow subgroups are $CP_2$-groups. In particular, a $p$-group $G$ is an $LCM$-group if and only if $G$ is a $CP_2$-group.}\\

\textbf{Lemma 2.3.} \textit{The class of $LCM$-groups is closed under taking subgroups and direct products.}\\

\textbf{Proof.} The fact that the $LCM$ property is inherited by the subgroups of an $LCM$-group is immediate. 

Let $G$ and $H$ be $LCM$-groups and $(g, h)\in G\times H$. It suffices to show that $(g, h)\in LCM(G\times H)$. We have
$$(g, h)\in LCM(G\times H)\Longleftrightarrow o((g, h)^n(g', h'))\mid \lcm(o((g, h)^n), o(g', h')), \forall \ n\in\mathbb{N}, \forall \ (g', h')\in G\times H$$ 
\begin{align}\label{r1}
\Longleftrightarrow \lcm(o(g^ng'), o(h^nh')) \mid \lcm(o(g^n), o(g'), o(h^n), o(h')),\forall \ n\in\mathbb{N}, \forall \ (g', h')\in G\times H.
\end{align}
Since $G$ and $H$ are $LCM$-groups, we have 
$$o(g^ng')\mid \lcm(o(g^n), o(g')), \forall \ n\in\mathbb{N}, \forall \ g'\in G\text{ and }o(h^nh')\mid \lcm(o(h^n), o(h')), \forall \ n\in\mathbb{N}, \forall \ h'\in H.$$
Then for any $n\in\mathbb{N}$ and $(g', h')\in G\times H$, the number $\lcm(o(g^n), o(g'), o(h^n), o(h'))$ is a multiple of both $o(g^ng')$ and $o(h^nh')$. It follows that (\ref{r1}) holds.
\hfill\rule{1,5mm}{1,5mm}\\ 

By Theorem 3.14 of \cite{19} and Theorem D of \cite{20}, it is known that all regular $p$-groups (in particular, all abelian $p$-groups) are $CP_2$-groups. Hence, by Lemma 2.2, they are also $LCM$-groups. Moreover, by Lemma 2.3, all abelian groups are $LCM$-groups. This information is useful for proving Theorem 2.4, which is the main result of this section. It improves Corollary 2.15 of \cite{3} that characterizes a minimal non-$LCM$-group $G$ under the hypotheses that $G$ is non-nilpotent and $2 \nmid |Fit(G)|$. For instance $A_4$ is a minimal non-$LCM$-group that cannot be detected via the respective  corollary since $Fit(A_4)\cong C_2^2$. Such cases are now covered by Theorem 2.4.\\

\textbf{Theorem 2.4.} \textit{Let $G$ be a  group. Then $G$ is a minimal non-$LCM$-group if and only if it is either a minimal non-$LCM$-$p$-group or a Schmidt group $P\rtimes Q$ such that $P$ is an $LCM$-group.}\\

\textbf{Proof.} Suppose that $G$ is a minimal non-$LCM$-group. Then $G$ is a non-$LCM$-group and all its proper subgroups are $LCM$-groups. These subgroups are also nilpotent by Lemma 2.2. We distinguish two cases based on the nilpotency of $G$.

\textit{Case 1.} If $G$ is nilpotent, then $G\cong P_1 \times P_2 \times\ldots \times P_k,$ where $k\geq 1$ and $P_1, P_2, \ldots, P_k$ are the Sylow subgroups of $G$, all of them being $LCM$-groups. If $k\geq 2$, then $G$ is an $LCM$-group by Lemma 2.3, a contradiction. Hence $k=1$ and $G$ is a minimal non-$LCM$-$p$-group.

\textit{Case 2.} If $G$ is non-nilpotent, then $G$ is a Schmidt group. By Lemma 2.1 \textit{i)}, it follows that $G=P\rtimes Q$. Since $P$ is a proper subgroup of $G$, we have that $P$ is an $LCM$-group.

Conversely, if $G$ is a minimal non-$LCM$-$p$-group, then the conclusion follows. If $G=P\rtimes Q$ is a Schmidt group, then $G$ is non-nilpotent. Lemma 2.2 implies that $G$ is not an $LCM$-group. To finish the proof, by Lemma 2.3, it suffices to show that all maximal subgroups of $G$ are $LCM$-groups. Let $H$ be a maximal subgroup of $G$. By making use of Lemma 2.1 \textit{v)}, we distinguish two cases.

\textit{Case 1.} $H\cong P\times \langle y^q \rangle$: $P$ is an $LCM$-group by our hypotheses, while $\langle y^q\rangle$ is also an $LCM$-group since it is a cyclic $q$-subgroup. By Lemma 2.3, we obtain that $H$ is an $LCM$-group.

\textit{Case 2.} $H\cong \Phi(P)Q_i$ with $i\in \{ 1, 2, \ldots, n_q \}$: It is clear that $\Phi(P)\cap Q_i=\{ 1\}$. In addition, by Lemma 2.1 \textit{ii)}, we deduce that all elements of $\Phi(P)$ are central in $G$, so $[g, h]=1$ for all $g\in \Phi(P)$ and $h\in Q_i$. It follows that $H\cong\Phi(P)\times Q_i$. We conclude that $H$ is an $LCM$-group, as a direct product of abelian $LCM$-groups.
\hfill\rule{1,5mm}{1,5mm}\\

Our previous result could be improved by providing information on the structure of a minimal non-$LCM$-$p$-group. This is left as an open question. Some comments related to it are written afterwards. These could be used as starting points for further research.\\

\textbf{Question 2.5.} \textit{Describe the structure of a minimal non-$LCM$-$p$-group.}\\

The following result shows that checking the $LCM$ property of a $p$-group is equivalent to verifying other properties that, in specific contexts, are easier to handle.\\

\textbf{Theorem 2.6} \textit{Let $G$ be a $p$-group. The following statements are equivalent:
\begin{itemize}
\item[i)] $G$ is an $LCM$-group;
\item[ii)] $G$ is a $CP_2$-group;
\item[iii)] $\Omega_k(G)=\{ x\in G \mid x^{p^k}=1\}, \forall \ k\in\mathbb{N}$;
\item[iv)] $G$ is a $P_2^*$-group. 
\end{itemize}}
\textbf{Proof.} The first 3 statements are equivalent according to Lemma 2.2 and Theorem D of \cite{20}. Let $k\in\mathbb{N}$. 

$ iii)\Longrightarrow iv):$ The hypothesis leads to $o(x)\leq p^k$ for all $x\in \Omega_k(G)$. Hence, $\Omega_k(G)\subseteq P_2^k(G)$. The other inclusion is trivial, so $G$ is a $P_2^*$-group.

$ iv)\Longrightarrow iii):$ We have that $P_2^k(G)=\Omega_k(G)$. It follows easily that $x^{p^k}=1$ for all $x\in\Omega_k(G)$, so $\Omega_k(G)\subseteq \{ x\in G \mid x^{p^k}=1\}$. The other inclusion is clear.
\hfill\rule{1,5mm}{1,5mm}\\

As a consequence of Theorem 2.6, one can say that the classes of minimal non-$LCM$-$p$-groups, minimal non-$CP_2$-$p$-groups and minimal non-$P_2^*$-$p$-groups coincide.\\

\textbf{Proposition 2.7.} \textit{Let $G$ be a non-$P_2$-$p$-group such that all its proper sections are $P_2$-groups. Then $G$ is a minimal non-$LCM$-$p$-group.}\\

\textbf{Proof.} Using the hypothesis, we infer that $G$ is a non-$P_2^*$-$p$-group, but all its proper sections are $P_2^*$-groups. By Theorem 2.6, $G$ is a non-$LCM$-$p$-group and all its proper sections (in particular, all its proper subgroups) are $LCM$-$p$-groups. The conclusion follows. 
\hfill\rule{1,5mm}{1,5mm}\\

It is worth mentioning that information on the structure of a non-$P_2$-$p$-group all of whose proper sections are $P_2$-groups is given by Theorem 6 of \cite{15}. Therefore, one may expect that some of those listed properties may be also satisfied by some minimal non-$LCM$-$p$-groups. The converse of Proposition 2.7 does not hold in general. For instance, $Q_{16}$ is the only minimal non-$LCM$-2-group of order $16$. Hence, it is also a non-$P_2$-$2$-group since $Q_{16}$ is not a $P_2^*$-group. However, $Q_{16}$ has a proper section isomorphic to $D_8$, which is not a $P_2^*$-group and therefore not a $P_2$-group. 

Let $G$ be a $p$-group of order $p^k$, with   $k\geq 1$. Note that if $k\leq p$, then $G$ is regular by 10.2 Satz \textit{b)} of \cite{12}. Hence, $G$ is an $LCM$-group. Based on this fact and taking into account Theorem 2.6, it is easy to show that the converse of Proposition 2.7 holds for $p$-groups of order $p^{p+1}$. 

By searching in GAP's \cite{22} SmallGroup library, we determined the number of minimal non-$LCM$-$p$-groups of order $p^k$ for specific values of $p$ and $k$. As we explained in the previous paragraph, it suffices to search for examples under the hypothesis $k>p$. Our findings are summarized below.
\begin{center}
\noindent\begin{tabular}{ |p{0.25cm}|p{0.2cm}|p{2.8cm}|p{9.9cm}|}
  \hline
  \multicolumn{4}{|c|}{$G$ is a minimal non-$LCM$-$p$-group of order $p^k$}\\
  \hline
$p$ & $k$ & no. of examples & $G\cong SmallGroup(p^k, i)$, with $i\in I$\\
 \hline
 2 & 3 & 1 & $I=\{3\}$\\
 \hline
 2 & 4 & 1 & $I=\{9\}$\\
\hline 
 2 & 5 & 3 & $I=\{10, 13, 14\}$\\
\hline 
 2 & 6 & 5 & $I=\{9, 18, 20, 21, 49 \}$\\
\hline 
 2 & 7 & 13 & $I=\{ 16, 21, 26, 29, 76, 83, 84, 100, 103, 104, 113, 117, 120\}$\\
 \hline
 2 & 8 & 48 & $I=\{23, 142, \ldots, 147, 160, 161, 164, 165, 201, 207, 208, 211, 212,$ $219, 222, 229, \ldots, 234, 238, 245, 249, 253, 264, 265, 268, 269, 276, \ldots ,$ $283, 285, 286, 288, 291, 414, 421, 450, 452\}$ \\
 \hline
 3 & 4 & 2 & $I=\{7, 9\}$\\
 \hline
 3 & 5 & 1 & $I=\{3\}$\\
\hline
 3 & 6 & 3 & $I=\{34, 96, 101 \}$\\
 \hline
 3 & 7 & 45 & $I=\{ 243, 244, 246, 250, 251, 252, 256, 257, 258, 274, 275, 280, 281,$ $282, 288, 289, 290, 292, \ldots, 299, 301, \ldots, 306, 350, 351, 352, 354, 357,$ $358, 359, 361, 365, 369, 370, 374, 375, 379\}$ \\
 \hline
 5 & 6 & 15 & $I=\{ 630, 631, 636, 639, 643, 645, 646, 647, 651, 652, 656, 661, 663,$ $664, 665\}$\\
 \hline
 \end{tabular}
 \end{center}

\section{Criteria for nilpotency and for the   $LCM$ property}

We need a series of preliminary results. 
By following the notations established in \cite{3}, given a group $G$, its subgroup generated by $LCM(G)$ is denoted by $LC(G)$. The first two results provide useful information on some properties of $LC(G)$ (see Theorem 2.10 of \cite{3}) and $LCM(G)$. The 3rd lemma indicates that, up to isomorphism, there is a unique non-$P_2$-2-group all of whose proper sections are $P_2$-groups (see the 1st paragraph of section 4 in \cite{15}). Finally, the 4th lemma characterizes the set $LCM(G\times H)$ under the hypothesis that $H$ is an $LCM$-group.\\  

\textbf{Lemma 3.1.} \textit{Let $G$ be a  group. Then $LC(G)$ is a nilpotent subgroup of $G$.}\\

\textbf{Lemma 3.2.} \textit{Let $G=K\rtimes H$ be a Frobenius group.
\begin{itemize}
\item[i)] If all Sylow subgroups of $K$ are $CP_2$-groups, then $LCM(G)=K$;
\item[ii)] If $K$ is abelian, then $LCM(G)=K$.
\end{itemize}}

\textbf{Proof.} Let $P_1, P_2, \ldots, P_k$ be the Sylow subgroups of $K$.

\textit{i)} By Example 2.3 of \cite{3}, since $P_1, P_2,\ldots, P_k$ are $CP_2$-groups, it follows that all these subgroups are included in $LCM(G)$. Let $x\in K$. Since $K$ is nilpotent, we may write $x=x_1x_2\ldots x_k$, where $x_i\in P_i, \forall \ i\in\{ 1,2,\ldots, k\}$. We have $x_1, x_2\in LCM(G)$. Also, since $(o(x_1), o(x_2))=1$ and $x_1x_2=x_2x_1$, it follows that $$o(x_1x_2)=o(x_1)o(x_2)=\lcm(o(x_1), o(x_2)).$$ By Remark 3.6 of \cite{3}, we conclude that $x_1x_2\in LCM(G)$. After repeating this argument for a finite number of times, we get $x\in LCM(G)$, so $K\subseteq LCM(G)$. 

Conversely, assume that $LCM(G)\not\subseteq K$. Let $x\in LCM(G)\setminus K$ and take $y=x^{-1}k$, where $k\in K$ and $k\neq 1$. It is easy to check that $y\in G\setminus K$. As both $x$ and $y$ do not belong to the Frobenius kernel $K$, they are conjugate to elements from $H$. Then $\lcm(o(x), o(y)) \mid |H|$. Since $x\in LCM(G)$, we have $o(xy)\mid\lcm(o(x), o(y))$. We conclude that $o(k)\mid |H|$. Since $o(k)\mid |K|$ also holds and $K$ and $H$ are of coprime orders, it follows that $o(k)=1$, a contradiction. Hence, $LCM(G)\subseteq K$.

\textit{ii)} Since $K$ is abelian, its Sylow subgroups are $CP_2$-groups. The conclusion follows based on \textit{i)}.
\hfill\rule{1,5mm}{1,5mm}\\

\textbf{Lemma 3.3.} \textit{$D_8$ is the only non-$P_2$-2-group all of whose proper sections are $P_2$-groups.}\\

\textbf{Lemma 3.4.} \textit{Let $G, H$ be groups. If $H$ is an $LCM$-group, then 
\begin{align}\label{r2}
(g, h)\in LCM(G\times H)\Longleftrightarrow o(g^ng') \mid \lcm(o(g^n), o(g'), o(h^n), o(h')),\forall \ n\in\mathbb{N}, \forall \ (g', h')\in G\times H.
\end{align}}
\textbf{Proof.} By using the hypothesis and the fact that $(g, h)\in LCM(G\times H)$ is equivalent to (\ref{r1}) (check the proof of Lemma 2.3), the conclusion follows.
\hfill\rule{1,5mm}{1,5mm}\\

The group $G$ that appears in the following result is one of the first two minimal non-$LCM$-3-groups (SmallGroup(81, 7)). We mention that $\expo(G)=9$.\\

\textbf{Lemma 3.5.} \textit{Let $k\geq 1$ be an integer and
$$G=\langle a,b,c, d \mid a^3=b^3=c^3=d^3=1, ba=abc, ca=acd, da=ad, cb=bc, db=bd, dc=cd  \rangle.$$  
\begin{itemize}
\item[i)] If $G_k=D_8\times C_{2^k}$, then $|LCM(G_k)|=|G_k|-8$;
\item[ii)] If $L_k=G\times C_{3^k}$, then $|LCM(L_k)|=|L_k|-108$.
\end{itemize}}
\textbf{Proof.} By using GAP, one can check that both results hold for $k=1$. Hence, we further assume that $k\geq 2$.

\textit{i)} Let $C_{2^k}=\langle t \rangle$. We use the notations established in the first section for $D_8$. Suppose that $n=1$. In this specific case, by using (\ref{r2}), it is easy to check that the following 4-tuples indicate the 8 elements $(g, h)\in G_k\setminus LCM(G_k)$: 
\begin{align*}
(g, h, g', h')\in &\{ (b, 1, a^3b, 1), (b, t^{2^{k-1}}, a^3b, 1), (a^2b, 1, a^3b, 1), (a^2b, t^{2^{k-1}}, a^3b, 1)\\ & (a^3b, 1, b, 1), (a^3b, t^{2^{k-1}}, b, 1), (ab, 1, b, 1), (ab, t^{2^{k-1}}, b, 1)\}.
\end{align*}

It remains to show that all other elements $(g, h)\in G_k$ satisfy (\ref{r2}). If $g\in \langle a\rangle$, then $g\in LCM(D_8)$. Then  (\ref{r2}) holds, so $(g, h)\in LCM(G_k)$ for all $h\in C_{2^k}$. If $g\not\in\langle a\rangle$, then $o(g)=2$. Also, $4\mid o(h)$ by the previous removal of 8 elements. If $n$ is odd, then $o(h^n)=o(h)$, and (\ref{r2}) becomes $$o(gg')\mid\lcm (2, o(g'), o(h), o(h')),\forall \ (g', h')\in G_k,$$ which holds since $o(gg')\mid 4$ and $4\mid\lcm (o(g), o(g'), o(h), o(h'))$. If $n$ is even, then (\ref{r2}) becomes
$$o(g')\mid \lcm(2, o(g'), o(h^n), o(h')), \forall \ (g', h')\in G_k,$$
which is also true. Hence, $|LCM(G_k)|=|G_k|-8$, as stated.   

\textit{ii)} We follow a similar argument as above. Let $C_{3^k}=\langle t\rangle$. There are 45 elements $g\in LCM(G)$. In these cases, we have $(g, h)\in LCM(L_k)$ for all $h\in C_{3^k}$. 

It remains to check the elements $(g, h)\in L_k$, with $g\in G\setminus LCM(G)=A\cup B$, where
$$A=\{ a, ad, ad^2, ac, acd, acd^2, ac^2, ac^2d, ac^2d^2, a^2, a^2d, a^2d^2, a^2c, 
a^2cd, a^2cd^2, a^2c^2, a^2c^2d, a^2c^2d^2 \},$$
$$B=\{b, bd, bd^2, bc, bcd, bcd^2, bc^2, bc^2d, bc^2d^2, b^2, b^2d, b^2d^2, b^2c, b^2cd, b^2cd^2, b^2c^2, b^2c^2d, b^2c^2d^2\}.$$ 
All of the above elements are of order 3. Due to similarity, we choose to treat only the case $g\in B$. Suppose that $n=1$. By  (\ref{r2}), the 4-tuples that indicate the 54 elements $(g, h)\in L_k\setminus LCM(L_k)$ are $(g, h, g', h')\in\{(\alpha, \beta, a, 1) \mid \alpha\in B, \beta\in \langle t^{3^{k-1}}\rangle\}$. Finally, we show that $$(g, h)\in LCM(L_k), \forall \ g\in B, h\not\in \langle t^{3^{k-1}}\rangle.$$ In this case, it is clear that $o(h)=3^s$ with $s\geq 2$. If $3\mid n$, then (\ref{r2}) becomes $$o(g')\mid \lcm(1, o(g'), 3^{sn}, o(h')), \forall \ (g', h')\in L_k,$$ which is true. Assume that $3\nmid n$. Then $o(g^n)=o(g), o(h^n)=o(h)$ and (\ref{r2}) is rewritten as $$o(g^ng')\mid \lcm(3, o(g'), 3^s, o(h')), \forall \ (g', h')\in L_k.$$
This also holds since $o(g^ng')\mid 9$ and $9\mid\lcm(3, o(g'), 3^s, o(h'))$.   

The other 54 elements $(g, h)\in L_k\setminus LCM(L_k)$ are identified by studying the case $g\in  A$.
\hfill\rule{1,5mm}{1,5mm}\\

To measure how close a group $G$ is to being an $LCM$-group, we introduce the following ratio:
$$\lcm(G)=\frac{|LCM(G)|}{|G|}.$$
Obviously, $\lcm(G)\in (0, 1]$ and $\lcm(G)=1$ if and only if $G$ is an $LCM$-group. If $G$ is a $p$-group, then $\lcm(G)=1$ if and only if any of the 4 equivalent statements of Theorem 2.6 holds. The following result shows that, by choosing any constant $c\in (0, 1)$, we cannot detect the $LCM$ property of a group $G$ of even/odd order by imposing a condition like $\lcm(G)>c$.\\

\textbf{Theorem 3.6.} \textit{There is no constant $c\in (0, 1)$ such that, whether we restrict to groups of even order or to groups of odd order, the condition  $\lcm(G)>c$  implies that $G$ is an $LCM$-group.}\\

\textbf{Proof.} By following the notations established in Lemma 3.5, there is a  sequence $(G_k)_{k\geq 1}$ of non-$LCM$-groups of even order and a sequence $(L_k)_{k\geq 1}$ of non-$LCM$-groups of odd order such that
$$\displaystyle\lim_{k\to\infty} \lcm(G_k)=\displaystyle\lim_{k\to\infty} \lcm(L_k)=1.$$
Hence, for any $c\in (0, 1)$ we can always find a non-$LCM$-group $G$ of even/odd order such that $\lcm(G)>c.$
\hfill\rule{1,5mm}{1,5mm}\\

However, for suitable choices of $c\in (0, 1)$, criteria for the $LCM$ property can be established by replacing the condition $\lcm(G)>c$ with $\lcm^*(G)>c$, where 
$$\lcm^*(G)=\min\{\lcm(S)\mid S \text{ is a section of } G\}.$$
The idea of working with sections is justified by the connections between the $LCM$ and $P_2^*$ properties in the case of $p$-groups, with the latter being used to define the concept of $P_2$-group, as discussed in the first two sections. Also, given a section $S$ of $G$, by using an isomorphism theorem, we infer that a section of $S$ is isomorphic to a section of $G$. This leads to
\begin{align}\label{r3}
\lcm^*(S)\geq \lcm^*(G) \text{ for any section } S \text{ of } G.
\end{align}

The following theorem is the main result of this section and it includes some criteria for nilpotency and for the $LCM$ property of a group. We will see that the nilpotency of a group can be recognized even without restricting to sections.\\

\textbf{Theorem 3.7.} \textit{Let $G$ be a group.
\begin{itemize}
\item[i)] If $\lcm(G)>\frac{1}{2}=\lcm(S_3)$, then $G$ is a nilpotent group;
\item[ii)] If $\lcm^*(G)>\frac{1}{2}=\lcm^*(S_3)$, then $G$ is a nilpotent group;
\item[iii)] If $G$ is a 2-group and $\lcm^*(G)>\frac{1}{2}=\lcm^*(D_8)$, then $G$ is an $LCM$-group;
\item[iv)] If $|\pi(G)|\geq 2$ and $\lcm^*(G)>\frac{1}{2}=\lcm^*(S_3)$, then $G$ is an $LCM$-group.
\end{itemize}} 

\textbf{Proof.} \textit{i)} We have 
$$\frac{|LC(G)|}{|G|}=\frac{|\langle LCM(G) \rangle|}{|G|}\geq\frac{|LCM(G)|}{|G|}=\lcm(G)>\frac{1}{2}.$$
Then $|LC(G)|>\frac{1}{2}|G|$. By using Lagrange's Theorem and Lemma 3.1, it follows that $G=LC(G)$ is a nilpotent group.

\textit{ii)} This is a straightforward consequence of \textit{i)}.

\textit{iii)} We argue by contradiction. Assume that $G$ is a non-$LCM$-group of minimal order with respect to the condition $\lcm^*(G)>\frac{1}{2}$. Let $S$ be a proper section of $G$. By (\ref{r3}), we have $$\lcm^*(S)\geq \lcm^*(G)>\frac{1}{2},$$ so $S$ is an $LCM$-group by minimality. According to Theorem 2.6, it follows that $G$ is a non-$P_2^*$-group and all its proper sections are $P_2^*$-groups. Then $G$ is a non-$P_2$-group and all its proper sections are $P_2$-groups.  By Lemma 3.3, we deduce that $G\cong D_8$, so $\lcm^*(G)=\frac{1}{2}$, a contradiction. 

\textit{iv)} Suppose that $G$ is a non-nilpotent group of minimal order with respect to the condition $\lcm^*(G)>\frac{1}{2}$. Since any proper subgroup $H$ of $G$ is isomorphic to a section of $G$, we make use of (\ref{r3}) to conclude that $H$ is nilpotent. Then $G$ is a Schmidt group. By following the notations established in Lemma 2.1, we have that $G=P\rtimes Q$ has a section
$S=G/Z(G)$ of order $p^rq$. Moreover, $S$ is a Frobenius group with kernel $K\cong C_p^r$ and complement $H\cong C_q$. By Lemma 3.2, we have that $LCM(S)=K$. Then $$\lcm^*(S)\leq\lcm(S)=\frac{|K|}{|K||H|}=\frac{1}{q}\leq \frac{1}{2}.$$
On the other hand, since $\lcm^*(G)>\frac{1}{2}$, we have $\lcm^*(S)>\frac{1}{2}$ as a consequence of (\ref{r3}). Thus we get a contradiction.    
\hfill\rule{1,5mm}{1,5mm}\\

For all previous criteria, the constant $\frac{1}{2}$ is the best possible one since $S_3$ is a non-nilpotent group (hence, a non-$LCM$-group by Lemma 2.2) and $D_8$ is a non-$LCM$-2-group. Also, criterion \textit{iii)} can be used to recognize the $CP_2$ and $P_2^*$ properties of any 2-group since they  are equivalent to the $LCM$ property by Theorem 2.6. 

Related to the previous theorem, we pose the following question that is briefly discussed afterwards.\\

\textbf{Question 3.8.} \textit{Let $p$ be an odd prime and let $G$ be a $p$-group. Is there a constant $c\in (0, 1)$ such that if $\lcm^*(G)>c$, then $G$ is an $LCM$-group?}\\

If $c$ exists, then $c\geq\frac{17}{25}$. This holds since $\frac{17}{25}$ is the maximum value among the values of $\lcm^*$ corresponding to the minimal non-$LCM$-$p$-groups that were indicated at the end of section 2. More exactly, for $G\cong SmallGroup(5^6, i)$ with $i\in\{631, 639, 643, 646, 647, 651, 656, 661, 664, 665 \}$, we have $\lcm^*(G)=\lcm(G)=\frac{17}{25}$.

\section{Minimum $\mathcal{F}$-covers}

Let $\mathcal{F}$ be a finite set of groups. In \cite{8}, the authors proved that if $\mathcal{F}$ is formed of nilpotent groups, then a nilpotent minimum $\mathcal{F}$-cover exists and there is a method to effectively construct it (see Theorem 7.3 of \cite{8} and its proof). A similar result holds if we replace ``nilpotent" by ``cyclic" in the previous phrase (see Theorem 7.1 of the same paper). However, if $\mathcal{F}$ is formed of solvable groups, then it is not always true that a minimum $\mathcal{F}$-cover is solvable. Indeed, an example given in the same paper shows that if $\mathcal{F}=\{ A_4, D_{10}\}$, then its unique minimum $\mathcal{F}$-cover is $A_5$. All these results and remarks are related to Question 8.2 of \cite{8}: 

\textit{For which classes $\mathcal{X}$ of groups, closed under taking   subgroups and direct products, is it true that, if $\mathcal{F}$ is a finite set of $\mathcal{X}$-groups, then there is a minimum $\mathcal{F}$-cover which is an $\mathcal{X}$-group?} 

At the time of writing this paper, this question is still unanswered when $\mathcal{X}$ is taken to be the class of abelian groups. In this section, our main objective is to show that if we choose $\mathcal{X}$ to be the class of $LCM$-groups (closed under taking subgroups and direct products by Lemma 2.3), then one finds  infinitely many sets $\mathcal{F}_k$ of $LCM$-groups such that there are exactly two  minimum $\mathcal{F}_k$-covers which are non-$LCM$-groups. Consequently, we point out a class of nilpotent groups that includes abelian groups and for which the answer to the above question is negative. 

We need some preliminary results. We use the notations established in the first section in their writing. The first two results outline the classification of non-abelian 2-groups possessing a cyclic maximal subgroup and provide information on the properties and the subgroups of certain groups of this class. These results  follow from Exercise 8a from section 1 of \cite{7}, Theorem 4.1 and section (4.2) (iii) of \cite{19}. The third lemma refers to $\mathcal{F}$-covers in the case of $p$-groups and is a consequence of Proposition 2.5 of \cite{8}.\\

\textbf{Lemma 4.1.} \textit{Let $G$ be a non-abelian 2-group having a cyclic maximal subgroup. Then $G$ is isomorphic to $M_{2^k}$, $D_{2^k}$, $Q_{2^k}$ or $QD_{2^k}$.}\\

\textbf{Lemma 4.2.} \textit{Let $k\geq 4$ be an integer. Then the following hold: 
\begin{itemize}
\item[i)] $M_{2^k}$ is a minimal non-abelian $2$-group;
\item[ii)] $D_{2^k}$ has only two non-cyclic maximal subgroups and both are isomorphic to $D_{2^{k-1}}$; 
\item[iii)] $Q_{2^k}$ has only two non-cyclic maximal subgroups and both are isomorphic to $Q_{2^{k-1}}$;
\item[iv)] $QD_{2^k}$ has only two non-cyclic maximal subgroups and they are isomorphic to $D_{2^{k-1}}$ and $Q_{2^{k-1}}$, respectively.
\end{itemize}
}

\textbf{Lemma 4.3.} \textit{Let $\mathcal{F}$ be a set of $p$-groups. Then every minimum $\mathcal{F}$-cover is a $p$-group.}\\

Next, we state and prove the main result of this section.\\

\textbf{Theorem 4.4.} \textit{Let $k\geq 4$ be an integer and let $\mathcal{F}_k=\{Q_8, C_{2^{k-1}} \}$ be a set of $LCM$-groups. Then  the only minimum $\mathcal{F}_k$-covers are the non-$LCM$-groups $Q_{2^k}$ and $QD_{2^k}$.}\\

\textbf{Proof.} Let $G$ be a minimum $\mathcal{F}_k$-cover. Then $G$ is a $2$-group according to Lemma 4.3. Since $G$ has a subgroup $F_2\cong C_{2^{k-1}}$, it follows that $2^{k-1}\mid |G|$. We may assume that $|G|=2^k$ since $G$ is non-cyclic, as it has a subgroup $F_1\cong Q_8$. It follows that $G$ is a 2-group possessing a cyclic maximal subgroup, so $G$ is isomorphic to one of the groups listed in Lemma 4.1. 

Suppose that $G\cong M_{2^k}$. Then $G$ is minimal non-abelian by Lemma 4.2. Hence, it has no subgroups isomorphic to $Q_8$, a contradiction. Let $I=\{3, 4, \ldots, k-1\}$. Suppose that $G\cong D_{2^k}$. By using Lemma 4.2, for each $i\in I$, the subgroups of order $2^i$ of $G$ are isomorphic to $C_{2^i}$ or $D_{2^i}$, so $G$ has no subgroups isomorphic to $Q_8$, a contradiction. It remains that $G\cong Q_{2^k}$ or $G\cong QD_{2^k}$. In both cases, by Lemma 4.2, we deduce that $G$ has subgroups isomorphic to $Q_{2^i}$ for all $i\in I$. Therefore, $G$ has subgroups isomorphic to $Q_8$. Hence, both $Q_{2^k}$ and $QD_{2^k}$ are confirmed to be minimum $\mathcal{F}_k$-covers, as stated. Moreover, by Lemma 2.3, both are non-$LCM$-groups since they have subgroups isomorphic to $Q_{16}$, which is the only minimal non-$LCM$-group of order 16.   
\hfill\rule{1,5mm}{1,5mm}\\

We end this section by noticing that Theorem 4.4 remains valid if ``$LCM$-groups" is replaced with ``$LCM$-2-groups". In this case, by Theorem 2.6, the result would also hold if the $LCM$ property is replaced with the $CP_2$ or $P_2^*$ property. A reduction of Theorem 4.4 to $p$-groups, with odd $p$, could be further studied. In this regard, we pose the following open problem:\\

\textbf{Question 4.5.} \textit{Let $p$ be an odd prime. Can one find infinitely many sets $\mathcal{F}$ of $LCM$-$p$-groups such that all their minimum $\mathcal{F}$-covers are non-$LCM$-$p$-groups?}\\

For instance, if $F$ is the non-abelian group of order $27$ and exponent 3, and $\mathcal{F}=\{F, C_9^2\}$, then its minimum $\mathcal{F}$-covers are two semidirect products of type $C_9^2\rtimes C_3$ (SmallGroup(243, 25) and SmallGroup(243, 26)). These details can be checked via GAP.

\section{Answering Amiri's questions}

Amiri's work \cite{2} includes the following main results which connect the $LCM$ and commutativity properties of groups of the same size with their sums of element orders and order sequences (see Theorems 3.8 and 3.9 of \cite{2}).\\

\textbf{Theorem 5.1.} \textit{Let $G$ and $H$ be $LCM$-groups of the same order. Then $\psi(G)=\psi(H)$ if and only if $\os(G)=\os(H)$.}\\

\textbf{Theorem 5.2.} \textit{Let $G$ and $H$ be abelian groups of the same order. Then the following statements are equivalent:
\begin{itemize}
\item[i)] The invariant factors of $G$ and $H$ are the same;
\item[ii)] $G\cong H;$
\item[iii)] $\os(G)=\os(H);$
\item[iv)] $\psi(G)=\psi(H).$
\end{itemize}}

Two open problems are formulated at the end of \cite{2}. In what follows, we recall these questions and provide an answer for each one.\\

\textbf{Question 5.3.} \textit{Let $G$ and $H$ be two groups of the same order such that $G$ is an $LCM$-group and $\expo(H)\mid\expo(G)$. If $\psi(G)=\psi(H)$, is it true that $\os(G)=\os(H)$?}\\  

By Theorem 5.1, the answer is affirmative if $H$ is an $LCM$-group. However, in general, the answer is negative. The first counterexamples are given by groups of order 256. More exactly, there are 1404 pairs $(G, H)$ of groups of order 256 satisfying the hypotheses of the previous question. The first 4 such pairs, with respect to the order of IDs in GAP's SmallGroup library, are given by the data below:
\begin{center}
\noindent\begin{tabular}{ |p{1.7cm}|p{0.95cm}|p{0.65cm}|p{2.5cm}|p{1.7cm}|p{0.95cm}|p{0.65cm}|p{2.5cm}|}
  \hline
IdGroup($G$) & $\expo(G)$ & $\psi(G)$ &$\os(G)$& IdGroup($H$) & $\expo(H)$ & $\psi(H)$ & $\os(H)$\\
 \hline
[256, 550] & 8 & 1519 & ((1, 1), (2, 7), (4, 120), (8, 128)) & [256, 2761] & 8 & 1519 & ((1, 1), (2, 71), (4, 24), (8, 160)) \\
\cline{5-5}
& & & & [256, 2768] & & &\\
\cline{5-5}
& & & & [256, 2790] & & &\\
\cline{5-5}
& & & & [256, 2799] & & &\\
\hline  
 \end{tabular}
 \end{center}
 
We mention that all 1404 pairs $(G, H)$  satisfy $\expo(G)=\expo(H)$, so $\expo(H)$ has never been a proper divisor of $\expo(G)$. Also, one can construct infinitely many counterexamples by choosing a pair $(G, H)$ given by the previous table and taking $(G\times C_p, H\times C_p)$ for any prime $p\geq 3$. 

Up to order 2000, we did not find pairs $(G, H)$ of groups of odd order satisfying all the hypotheses of Question 5.3 such that $\os(G)\neq\os(H)$. Hence, it would be an idea to study this question by adding this hypothesis which would imply that $G$ and $H$ are solvable groups.\\

\textbf{Question 5.4.} \textit{Let $G$ and $H$ be two groups of the same order such that $G$ is an $LCM$-group and $\expo(H)\mid\expo(G)$. If $H$ is a non-$LCM$-group, is it true that $\psi(H)<\psi(G)$?}\\

Again, in general, the answer is negative. Up to order 100, there are 515 counterexamples $(G, H)$. For 303 of them, we actually have $\psi(G)<\psi(H)$. The first 4 out of these 515 counterexamples are given by the data below:
\begin{center}
\noindent\begin{tabular}{ |p{1.7cm}|p{1cm}|p{1cm}|p{1.7cm}|p{1cm}|p{1cm}|}
  \hline
IdGroup($G$) & $\expo(G)$ & $\psi(G)$ & IdGroup($H$) & $\expo(H)$ & $\psi(H)$\\
 \hline
[16, 3] & 4 & 47 & [16, 13] & 4 & 47\\
\cline{1-1}
[16, 10] & & & & & \\
\hline
[32, 5] & 8 & 175 & [32, 38] & 8 & 175 \\
\hline
[32, 22] & 4 & 95 & [32, 6] & 4 & 103\\
\hline  
 \end{tabular}
 \end{center}
 
Assuming that $G$ and $H$ are of odd order, the first counterexamples are given by groups of order 243. There are 60 such instances $(G, H)$ and for 30 of them we have $\psi(G)<\psi(H)$. Here are the first 4 out of the 60 counterexamples:
\begin{center}
\noindent\begin{tabular}{ |p{1.7cm}|p{1cm}|p{1cm}|p{1.7cm}|p{1cm}|p{1cm}|}
  \hline
IdGroup($G$) & $\expo(G)$ & $\psi(G)$ & IdGroup($H$) & $\expo(H)$ & $\psi(H)$\\
 \hline
[243, 4] & 9 & 1699  & [243, 25] & 9 & 1807 \\
\cline{4-4}
 & &  & [243, 30] & & \\
\cline{4-4}\cline{6-6}
 & & & [243, 55] & & 1699  \\
\cline{4-4}
 & & & [243, 60] & & \\
\hline  
 \end{tabular}
 \end{center}

For all 575 counterexamples $(G, H)$, we have $\expo(G)=\expo(H)$ so, again, $\expo(H)$ has never been a proper divisor of $\expo(G)$. If $(G, H)$ is one of the above counterexamples, then $(G\times C_p, H\times C_p)$ is also a counterexample for any prime  $p\geq 5$.

\bigskip\noindent {\bf Acknowledgements.} The author is grateful to the reviewer for his/her remarks which improved the initial version of the paper. 
    
\bigskip\noindent {\bf Declarations}

\bigskip\noindent {\bf  Funding.} The author did not receive support from any organization for the submitted work.

\bigskip\noindent {\bf Conflicts of  interests.} The author declares that there is  no conflict of interest.

\bigskip\noindent {\bf  Data availability statement.} The manuscript has no associated data.

\vspace*{3ex}
\small
\hfill
\begin{minipage}[t]{7cm}
Mihai-Silviu Lazorec \\
Faculty of  Mathematics \\
"Al. I. Cuza" University \\
Ia\c si, Romania \\
e-mail: {\tt silviu.lazorec@uaic.ro}
\end{minipage}
\end{document}